\documentclass[openbib,leqno]{article}
\usepackage{graphicx,amsfonts}
\begin{document}

\title{Endomorphic Elements in Banach Algebras}
\author{V.A. Babalola  \\ Department of Mathematics, \\ University of Ibadan, \\
 Ibadan, Nigeria. \\
 va.babalola@mail.ui.edu.ng}
\date{}

\maketitle

\begin{abstract}
 The use of the properties of actions on an algebra to enrich the study of the algebra is
 well-trodden and still fashionable. Here, the notion and study of endomorphic elements of
 (Banach) algebras are introduced. This study is initiated, in the hope that it will open up,
 further, the structure of (Banach) algebras in general, enrich the study of
 endomorphisms and provide examples.

 In particular, here, we use it  to classify algebras for the
 convenience of our study.  We also present results on the
 structure of some classes of endomorphic elements and
 bring out the contrast with idempotents.
\end{abstract}

\noindent{\em 2000 MSC-class:} 08A35, 46H20 (primary); 47N99,
16N40 (secondary).

\section{Motivation}
\label{intro} Dictated by convenience, we shall introduce a twist
in the
    usual terminologies concerning elements. Let $A$ be a complex
    or real algebra. $a \in A$ will be called a left element, a
    right element or just an element according as $a$ is,
    considered as a self map of $A$, of the form $x \mapsto ax$, $x \mapsto xa$,
    or $x\mapsto axa$. A bi-element is a left and right element. A
    self map $T$ of $A$ is said to be endomorphic if it is linear
    and $T(ab) = T(a)T(b) \;\;\forall a, b \in A$. The meaning of terms like
    endomorphic left element, compact endomorphic right element, etc becomes immediately
    obvious.

    The disc algebra is a semisimple commutative Banach algebra having its maximal ideal space
 homeomorphic with  (0,1). It has no nontrivial endomorphic map. In 1980 Kamowitz
 \cite{kamowitz-80} posed the problem of studying compact endomorphic maps of Banach
 algebras in  general to see if the situation generalizes. In the above paper and again in
 1989 \cite{kam-sch-wort-89}, 1998 \cite{kamowitz-98}, 2000 \cite{fein-kamo-00}, and 2004
 \cite{fein-kamo-04}, along with others, he went on to discuss this for commutative
 semisimple algebras. To complement their work, I have relaxed the restriction on the
 algebra to accommodate all algebras and considered endomorphic elements in place
 of compact endomorphic maps.  This is with a view of throwing more light on the study of
 endomorphic maps in general.

  Another motivation for this attempt is a paper of Zemanek
\cite{zemanek-79} on idempotents which he presented to the
International Congress of Mathematicians in Helzinki in 1978.
There he discussed the
  structure of idempotents in Banach algebras in general and semisimple Banach algebras in
  particular. That endomorphic elements and idempotents are closely related is emphasized and
  exploited in our presentation.

 More particularly, this paper deals with the structure of $L(A)$, the class of
 endomorphic left elements in $A$ given by
 \[ L(A) := \{ a \in A: axay = axy \;\;\;\;\;\;\forall x,y \in A \}. \]
  Other possibility is the class $R(A)$ of all
 endomorphic right elements of $A$. This class is by implication
 already being treated as appropriate adjustments of the
 $L(A)$-situation will cover this.

 Unexplained terms are those of \cite{bon-dun-73}. Except
 otherwise explained, $A$ is an algebra.

 \section{Introduction of Appropriate Terms and
  Notations}
  \label{intro-terms-nota}
 \paragraph{Definition}  Let $A$ be an algebra. Then $I(A)$
  denotes the class of all idempotents of $A$.

\paragraph{Definition} ~ \\
 (1) {\sl Algebra Without Order:\ } An algebra $A$ is
  without order iff
    \[ (a \in A) \wedge [(ax = 0 \;\forall x \in A)\vee (xa = 0 \;\forall x \in
  A)] \Rightarrow  (a = 0). \]
 (2) {\sl Nice Algebra:\ } An algebra $A$ is
  nice iff
  \[ (a \in A) \wedge (axy = axay \;\;\forall x,y \in A) \\
 \Rightarrow (ax = axa \;\;\forall x \in A). \]
 (3) {\sl Very Nice Algebra:\ } An
 algebra $A$ is  very nice iff
 \[ (a \in A)\wedge(axy = axay \;\;\forall x,y \in A) \Rightarrow (a \in I(A)) \wedge
 (ax = axa \;\;\forall x \in A). \]


\[ 
\unitlength 1.00mm 
\linethickness{0.4pt}
\ifx\plotpoint\undefined\newsavebox{\plotpoint}\fi 
\begin{picture}(145.75,55.00)
\put(75.5,26.75){\line(0,1){5.5}} 
\put(73.5,26.75){\line(0,1){6.00}} 
\multiput(77,32.00)(-.0333333,.04){75}{\line(0,1){.04}} 
\multiput(72.25,32.50)(.0373134,.0335821){67}{\line(1,0){.0373134}} 
\put(104.5,26.75){\line(0,1){6.00}} 
\put(102.5,26.75){\line(0,1){6.00}} 
\multiput(106,32.00)(-.0333333,.04){75}{\line(0,1){.04}} 
\multiput(101.25,32.50)(.0373134,.0335821){67}{\line(1,0){.0373134}} 
\multiput(34.75,45.00)(-.0410448,-.0335821){67}{\line(-1,0){0.5}} 
\multiput(34.75,45.00)(-.05,.0333333){60}{\line(-1,0){.05}} 
\put(26.50,46){\line(1,0){5.5}} 
\put(26.50,44.00){\line(1,0){5.5}} 
\multiput(69.00,45.00)(-.0410448,-.0335821){67}{\line(-1,0){0.5}} 
\multiput(69.00,45.00)(-.05,.0333333){60}{\line(-1,0){.05}} 
\put(60,46){\line(1,0){6.0}} 
\put(60,44.00){\line(1,0){6.0}} 
\multiput(95.00,45.00)(-.0410448,-.0335821){67}{\line(-1,0){0.5}} 
\multiput(95.00,45.00)(-.05,.0333333){60}{\line(-1,0){.05}} 
\put(86,46){\line(1,0){6.0}} 
\put(86,44.00){\line(1,0){6.0}} 
\multiput(121.00,45.00)(-.0410448,-.0335821){67}{\line(-1,0){0.5}} 
\multiput(121.00,45.00)(-.05,.0333333){60}{\line(-1,0){.05}} 
\put(112,46){\line(1,0){6.0}} 
\put(112,44.00){\line(1,0){6.0}} 
\multiput(95.00,19.00)(-.0410448,-.0335821){67}{\line(-1,0){0.5}} 
\multiput(95.00,19.00)(-.05,.0333333){60}{\line(-1,0){.05}} 
\put(86,20){\line(1,0){6.0}} 
\put(86,18.00){\line(1,0){6.0}} 
\put(2,40){\framebox(24,10.75)[cc]{{\bf B* Algebra}}}
\put(35,34.5){\framebox(24.5,18.25)[cc]{\parbox{22.5mm}{\bf Algebra with \\
approximate \\ identity}}}
\put(69,35.5){\framebox(16.5,15.25)[cc]{\parbox{14mm}{\bf Very \\
nice \\ algebra}}}
\put(95.25,35.5){\framebox(15.75,15.5)[cc]{\parbox{13.5mm}{\bf
Algebra \\ without \\ order}}}
\put(121.25,39.25){\framebox(16.5,11.75)[cc]{\parbox{14mm}{\bf
Nice \\ algebra}}}
\put(95,14){\framebox(22.75,12.25)[cc]{\parbox{20mm}{\bf Semiprime
\\ algebra}}}
\put(61.5,13.25){\framebox(23.75,13)[cc]{\parbox{21.5mm}{\bf
Semisimple algebra}}} \put(67.5,6.5){\makebox(0,0)[cc] {{\large
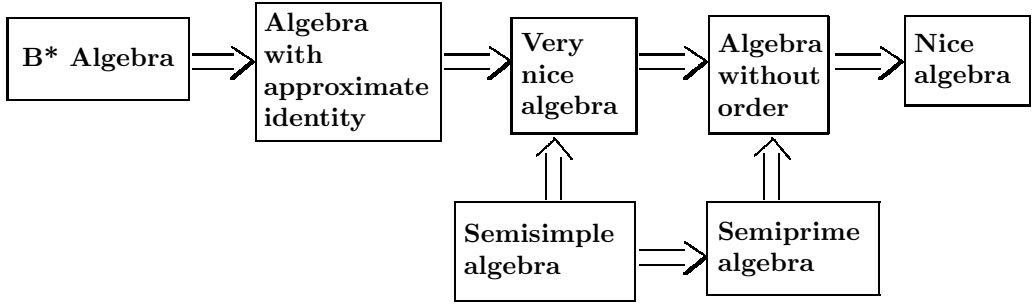
Figure 1: Relations Among Relevant Types of Algebras}}}
\end{picture} \]

  \paragraph{Definition} ~ \\
  (1) $N'_3(A)$ denotes the maximal
 nilpotent subalgebra of the algebra $A$ of order $3$. i.e.
 \[ N'_3(A) := \{a \in A : axy = 0 \;\;\forall x,y \in A\}. \]

 \noindent (2) $\displaystyle N_3(A):= \{a \in A : a^3 = 0\}$.

 \noindent (3) $\displaystyle N(A):= \{a \in A :  a^n = 0 \mbox{\ for some\ } n \in \mathbb{N} :=
 \{1, 2, 3, \ldots, n, \ldots\}\}$.

 \noindent (4) If $A$ is a normed algebra, then $QN(A)$ denotes the class of all
 quasinilpotent
 elements of $A$.

 \paragraph{Definition}
 $A$ is  an endomorphic left algebra iff $A = L(A)$.

 \paragraph{Definition}
  If $A$ is a normed algebra and $U$, $V \in A$
  then
  \[ d(U,V) := \inf_{x \in U,\; y\in V}\|x-y\|. \]


  \setcounter{paragraph}{0}
   \section{General Properties of Endomorphic Left Elements}
  \paragraph{Theorem} The following statements hold: \\
 (1) $(a \in L(A)) \Rightarrow$ $(\exists n \in \mathbb{N} \mbox{\ \ such that\ \ }1
 \leq n \leq3 \mbox{\ \ and\ \ } a^n = a^{n+1})$. \\[0.25em]
 (2) $L(A)I(A) \subset I(A)$. \\[0.25em]
 (3) $[L(A)]^2 \subset L(A)$. \\[0.25em]
 (4) $\displaystyle L(A) = \bigcup_{b \in G(\tilde{A})}b L(A)b^{-1}$ \\
 where $G(\tilde{A})$ is the class of all invertible elements of $\tilde{A}$, the minimal
 algebra which has identity and contains $A$. \\[0.25em]
  (5) If $A$ is a normed algebra then $L(A)$ is closed in $A$.

\noindent {\bf Proof\ \ } (1) -- (4) follow easily from the
definition of a endomorphic left element and manipulation of
standard techniques. Take, for example, (3) and (5).

Consider (3). Let $a, b \in L(A)$. Then
\[ abxaby = abxay = abxy \;\;\;\;\;\forall x,y \in A. \]
Therefore $ab \in L(A)$ and (3) holds.

That $L(A)$ is closed in $A$ follows from the joint continuity of
product in $A$. ~ \hfill $\Box$

\paragraph{Note} ~

 \noindent
 (1) If $(a \in L(A))\wedge(a \neq 0)\wedge (a \neq
1)$ then $a$ is a divisor of zero. \\[0.25em]
 (2) An immediate conclusion from (1) is the following: An algebra
without a divisor of zero has no non-trivial endomorphic elements. \\[0.25em]
 (3) The inclusion in 3.1(3) can be strict. Take a non-trivial
nilpotent algebra $A$ of order 2. Then $A = L(A)$ and $[L(A)]^2 =
\{0\}$. \\[0.25em]
 (4) If $A$ is nice then

 (i) $(a \in L(A)) \Rightarrow$ $(\exists n \in \mathbb{N} \mbox{\ \ such that\ \ }1
 \leq n \leq 2 \mbox{\ \ and\ \ } a^n = a^{n+1})$.

 (ii) $[L(A)]^2 = L(A) \cap I(A)$. \\[0.25em]
 {\textsl{Reason}.\ \ } (i) follows from the definition of a nice algebra.
 Consider (ii). From
 3.1(3) and the definition of a nice algebra, $[L(A)]^2 \subset L(A) \cap I(A)$. Take
 $a \in I(A)\cap L(A)$. Then $a^2 = a$. Therefore $a \in [L(A)]^2$. Therefore
 $[L(A)]^2 \supset L(A) \cap I(A)$. \\[0.25em]
 (5) Part of Zemanek's characterization in [7] of central idempotents is that for a
 semisimple Banach algebra $A$
 \[ (e \in Z(A)\cap I(A)) \Leftrightarrow (eI(A) \subset I(A)). \]
 There He also gave an example of a non-semisimple Banach algebra
 having a non-central idempotent $e$ with
 \[ eI(A) \subset I(A). \]
 In fact his $A$ is the Banach algebra of complex upper triangular
 $2 \times 2$ matrices $\left( \begin{array}{cc}
 a & b \\
 c & 0
 \end{array} \right)$ and his $e$ is
 $\left( \begin{array}{cc}
 1 & 0 \\
 0 & 0
 \end{array} \right)$.
 Indeed, this observation is part of a greater truth which is
 3.1(2). \\[0.25em]
(6) Equality holds in 3.1(2) if $I(A) \subset Z(A)$.

\noindent {\textsl{Reason}.\ \ } In this case $I(A) \subset L(A)$
and $a \in I(A) \subset L(A) \Rightarrow a =aa \in L(A) I(A)$.
Thus \[ L(A) I(A) \supset I(A). \]

\noindent (7) Interpreting  (6) gives the following: If in an
algebra $A$ every idempotent is in the centre of $A$
then every endomorphic left element can be factored as the
product, $ai$, of an endomorphic left element and an idempotent.

\setcounter{paragraph}{0}
\section{Nilpotency and Endomorphic Left Elements}
\paragraph{Theorem} The following statements are e\-qui\-valent:
\\
(1) $a \in  N(A) \bigcap L(A)$. \\[0.25em]
 (2) $a \in N_3(A)\bigcap L(A)$. \\[0.25em]
 (3)] $a \in  N'_3(A)$. \\[0.25em]
If $A$ is a normed algebra, then each of (1) -- (3) is equivalent
to the following: \\[0.25em]
(4) $a \in QN(A) \cap L(A)$.

 \noindent {\bf Proof\ \ }
  \begin{eqnarray*}
  \lefteqn{Q N(A) \cap L(A)  \subset N(A) \cap L(A)  \subset  N_3(A)\cap
  L(A)}
  \\
  {} &  \subset & N'_3(A)\subset N(A) \cap L(A) \subset Q N(A) \cap L(A)
  \end{eqnarray*}
  is immediate from definition
 and 3.1(1). The first and last terms are included only if $A$ is
 a normed algebra. ~ \hfill $\Box$

\paragraph{Theorem} $A$ is an endomorphic left algebra  iff $A =
N'_3(A)$. i.e. $A^3 = \{0\}$.

\noindent {\bf Proof\ \ } Let $K$ be the scalar field of $A$. Let
$\alpha \in K$ be arbitrary.
 Then
 \begin{eqnarray*}
 \lefteqn{(A= L(A))\wedge (\exists b\in A \mbox{\ with \ } b^3 \neq
 0)} \\
 & \Rightarrow & (b^4 = b^3) \wedge ((\alpha b)^4 = (\alpha b^3)), \;\;\;\mbox{\ from
 3.1} \\
 & \Rightarrow & (\|b^4\| = \|b^3\|) \wedge (|\alpha|^4 \|b^4\| = |\alpha|^3 \|b^3\|)
 \\
 & \Rightarrow & \alpha = \pm 1 \mbox{\ or\ } 0.
 \end{eqnarray*}
 This is a contradiction since $\alpha \in K$ is arbitrary. Hence
 \[ (L(A) = A) \Rightarrow (a^3 = 0 \;\; \forall a \in A). \]
 The reverse implication  follows from 4.1. ~ \hfill $\Box$

  \paragraph{Theorem}
  The following statements are e\-qui\-valent:

 \noindent
 (i) $N'_3(A) = \{0\}$. \\[0.25em]
 (ii) $A$ is without order.

 \noindent
 {\bf Proof\ \ }
 It is easy to see. ~ \hfill $\Box$

 \paragraph{Corollary} Let $A$ be an algebra without order and an endomorphic left algebra.
 Then
 $A$ is trivial.

 \paragraph{Note} If $A$ in Theorems 4.1 -- 4.3 is in addition nice then the suffix and
 superfix 3 can be replaced by 2 everywhere they appear in these theorems.
 Also in the sequel, all general results involving suffix 3 / superfix 3 are also true
 when the algebra is, in addition, assumed nice and the said 3 is replaced by 2.  All these
 can be seen with the help of 3.2(4)(i), the observation that for a nice algebra $A$,
 \[ N'_3(A) = N'_2(A),\;\;\;\; N_3(A) = N_2(A) \]
 and some manipulation. \hfill $\Box$


 \setcounter{paragraph}{0}
 \section{Connectedness of $L(A)$, A, a Banach Algebra}
 Throughout this section, $A$ is, among other things, a Banach
 algebra. We continue with our structure theorems under this
 additional assumption.
 \paragraph{Theorem}
  Let $A$, be a Banach algebra.  Then the following hold: \\[0.25em]
  (1) $N'_3(A)$ is the component of the origin in $L(A)$. \\[0.25em]
  (2) $N'_3(A)$ is unbounded if it is not a singleton. \\[0.25em]
  (3) $L(A)$ is disconnected iff
  \[ Q(A) := L(A)\setminus N'_3(A)  \neq \emptyset. \]
    Moreover, \\
    (4) ~ \hfill \parbox{100mm}{\[ Q(A) \neq \emptyset \Rightarrow
    d(N'_3(A), Q(A)) \geq 1. \]} \hfill ~

\noindent{\bf Proof\ \ }
    Consider (1). $N'_3(A) \subset C_0$, the component of $L(A)$
    containing the origin. This is because \\
    \[ (0 \neq a \in N'_3(A)) \Rightarrow (f := t \mapsto ta : [0,1]  \rightarrow N'_3(A)). \]
    Thereby joining the origin to $a$ through $N'_3(A)$. To establish
    $C_0 = N'_3(A)$ it remains to prove (i).  Suppose $a \in N'_3(A)$
    and $b \in Q(A)$. Then
    $(b - a)b^3 =b^3$ and so $\|b -a\| \geq 1$. Therefore
    $d(N'_3(A), Q(A)) \geq
    1$. Hence (1) holds.

    Consider (2). If $N'_3(A)$ is not a singleton then $\exists a
    \in N'_3(A)$ such that $a \neq 0$. Then $\{ta : t \in [0,
    \infty)\}$ is an unbounded subset of $N'_3(A)$, thereby making $N'_3(A)$, itself,
    unbounded.

     (3) will be established , in view of (4), if we show that:
     \begin{quote}
     ($L(A)$ is disconnected) $\Rightarrow$ ($L(A)$ has a
     non-quasinilpotent member).
     \end{quote}
     Now its contrapositive is:
     \[ (L(A) = N'_3(A)) \Rightarrow (L(A) \mbox{\ is connected})
     \]
     which is true by (1). ~ \hfill $\Box$ \\


    Our next main task is to characterize the isolated points of
    $L(A)$. To facilitate this we establish the following lemma.

    \paragraph{Lemma}
    Suppose $a,b \in L(A)$ with $ab = ba$, $b^3 \neq a$ and $a \in
    I(A)$. Then $\|a-b\| \geq 1$.

\noindent{\bf Proof\ \ }
    Given $a$ and $b$ satisfying the hypothesis of the lemma, then
      \[ \begin{array}{lcrcl}
     \mbox{\hspace{10.5ex}} 0 & \neq & (a - b)^3 & = & (a - b)^5. \\
    \mbox{Therefore\ }
     0 & \neq & \|(a - b)^3\| & \leq & \|(a -
    b)^2\|\|(a - b)^3\|. \\
    \mbox{Therefore}  & & \|a - b\| & \geq & 1.
    \end{array}  \]
    ~ \hfill $\Box$
    \paragraph{Theorem (Isolation)}
    The follow\-ing sta\-te\-ments are equi\-valent:
    \begin{enumerate}
    \item [(1)] $a$ is isolated in $L(A)$.
    \item [(2)] $a \in I(A)\cap Z(A)$ and any $b \in L(A)$ with
    $b^3 = a$ is such that $b=a$.
    \end{enumerate}

   \noindent{\bf Proof\ \ } Let $a \in L(A)$. Then the following conditions are mutually
    exclusive and exhaustive: \\[0.25em]
    (i) $a^2 \neq a$. \\[0.25em]
    (ii) $a^2 = a$ and $\exists \; b \in L(A)$ such that
    $a = b^3 =b^2 \neq b$. \\[0.25em]
    (iii) $a^2 = a$ and $\exists \; b \in L(A)$ such that
    $a = b^3 \neq b^2 = b$. \\[0.25em]
    (iv) $a^2 = a$ and every $b \in L(A)$ with
    $a = b^3$ is such that  $b^3 =b^2 = b$. \\[0.25em]
    For $a$ to be isolated in $L(A)$; (i), (ii) and (iii) are not
    possible as then $a$ will be arc connected in (i) to $a^2$ by
    the line segment $f_1\left|_{[0,1]} \right.$ of the ray
    \[ f_1 := t \mapsto ta + (1-t)a^2 : [0, \infty) \rightarrow
    L(A), \]
    in (ii)  and in (iii) to $b$ by the line segment
    $f_2\left|_{[0, 1]}\right.$ of the ray
    \[ f_2 := t \mapsto ta + (1-t)b : [0, \infty ) \rightarrow
    L(A). \]
    We are left with alternative (iv). Thus the conclusion here is
    that if $a$ is isolated in $L(A)$ then \\[0.25em]
    (iv)$'$ $a \in I(A)$ and  any $b \in L(A)$ with $b^3 =a$ is such that
    $b=a$. \\[0.25em]
    \indent Again let $a \in L(A)$. Then
    \begin{eqnarray*}
    F_a(\cdot, \cdot) \!&:\!= \!& \!(x, y) \mapsto a - xay - axay : A \!\times \!A
    \rightarrow A \\
    F'_a(\cdot, \cdot) \!& \!:= \!& \!(x, y) \mapsto a - xay - xaya : A \!\times \!A
    \rightarrow A
    \end{eqnarray*}
    are both continuous maps with $F_a(A\times A) \subset L(A)$,
     $F'_a(A\times A) \subset L(A)$ and
     \[ F_a(0, 0) = a = F'_a(0, 0). \]
     If, in addition, $a$ is isolated in $L(A)$, it follows that \\[0.25em]
     (v) $axay = xay \;\;\;\;\;\forall x,y \in A$; \\[0.25em]
     (vi) $xaya = xay \;\;\;\;\;\forall x,y \in A$.

     Now manipulation of (iv)$'$, (v) and (vi) shows that:
     \begin{quote}
      ($a$ is
     isolated in $L(A)$) $\Rightarrow$ $(a \in Z(A))$.
     \end{quote}

     Altogether therefore  $(1) \Rightarrow (2)$.

     Conversely assume (2). Note that in (2)
     \[ [(b^3 = a) \Rightarrow (b = a)] \Leftrightarrow [(b \neq a) \Rightarrow (b^3 \neq
     a)]. \]
     Therefore along with Lemma 5.2,
     \[ (2) \Rightarrow (b \in L(A) \mbox{\ with\ } b\neq a \mbox{\ is such that\ }
      \|b - a\| \geq     1). \]
     Hence
     $a$ satisfying (2) is isolated in $L(A)$. ~ \hfill $\Box$

 \paragraph{Theorem}
 A component of $L(A)$ is either unbounded or a singleton.

 \noindent{\bf Proof\ \ }
 Theorem 5.3 characterizes the singleton components of $L(A)$. For
 $a \in L(A)$, let $C_a$ be the component containing $a$.
 The conditions (i) -- (iii) in the proof of Theorem 5.3 are mutually  exclusive
 and exhaust the possibilities for $C_a$ to be a non-singleton component.

 \noindent Condition(i):\ \ With $f_1$ as defined in the proof of
 the the theorem, the entire ray $f_1([0, \infty))$ lie in
 $C_a$. Hence $C_a$ is unbounded.

 \noindent Conditions (ii) and (iii):\ \  Similarly the entire ray
 $f_2([0, \infty))$ lie in $C_a$ under each of the conditions (ii)
 and (iii), where $f_2$ is also as defined in the proof of Theorem
 5.3. Therefore $C_a$ is also unbounded under each of the conditions
 (ii) and (iii). ~ \hfill $\Box$

 \paragraph{Theorem}
 A non-central or non-idem\-po\-tent member of $L(A)$ is in some
 unbounded component of $L(A)$.

 \noindent{\bf Proof\ \ }
 The proof follows from the proof of Theorem 5.4. ~ \hfill  $\Box$

\paragraph{Theorem}
The component of $L(A)$ containing the origin is a singleton iff
$A$ has no order. Otherwise it is unbounded.

\noindent{\bf Proof\ \ } The theorem follows from Theorem 4.3 and
Theorem 5.1(2). ~ \hfill $\Box$

 \paragraph{Theorem (Existence of Isolated Points)} $L(A)$ has an iso\-la\-ted point iff $A$
 has no order.

 \noindent{\bf Proof\ \ }
  Consider the {\em only if} part. Suppose
 $a \in L(A)$ is isolated. Then $a \in I(A)\cap L(A)\cap
 Z(A)$ and every $b \in L(A)$ with $b^2 = a^2$ is such that $b =
 a$ Now $A$ is with order. Therefore $\exists$ $c \in A$ such that
 ($cx = xc = 0 \;\; \forall x \in A)\wedge (c \neq 0)$. Therefore $c \in
 L(A)$. Moreover $a + c \in L(A)$ since
 \begin{eqnarray*} (a +c)x(a + c)y & =  & (axa +axc + cxa + cxc)y \\
 & = & (ax + 0 + 0 + 0)y \\
 & = & axy + cxy \\
 & = & (a + c)xy \;\;\;\;\forall x, y \in A.
 \end{eqnarray*}
 Now $(a + c)^2 = a^2$ but $a + c \neq a$ since $c \neq 0$. This
 is a contradiction. Therefore $A$ is without order. As such $L(A)$ has an isolated
 point only if $A$ has no order.

 The {\em if} part follows from Theorem 4.3. ~ \hfill $\Box$

\paragraph{Corollary} If $A$ is a very nice Banach algebra,
then $a\in L(A)$ is isolated iff $a$ is central.


\setcounter{paragraph}{0}
 \section{Connectedness and  Endomorphic Left Elements; $A$, a Very Nice Banach Algebra}
 \paragraph{Theorem} Let $A$ be a very nice Banach algebra. Then the component of $L(A)$
 containing an element $a \in L(A)$ coincide with the component of $I(A)$ containing $a$.

\noindent{\bf Proof\ \ }
 Let $a \in L(A)$. Let $K_a$ be the
 component of $a$ in $I(A)$ and $C_a$ the component of $a$ in
 $L(A)$. Since $a \in L(A)$ and $K_a = \{\omega a {\omega}^{-1} :
 \omega \in G(\tilde{A}) \}$, it follows from Theorem 3.1(4) that $K_a
 \subset L(A)$. Therefore $K_a \subset C_a$. Since $L(A) \subset
 I(A)$ by definition, we conclude that $K_a = C_a$. ~ \hfill
 $\Box$ \\

    Let $A$ be a very nice Banach algebra. Having identified the
components of $L(A)$ as those of $I(A)$, the following conclusions
then follow from Zemanek's paper, \cite{zemanek-79}:
\begin{enumerate}
\item[(1)] $(e$, $f \in L(A)) \wedge (r(e - f) < 1) \Rightarrow$
($e$ and $f$ are arc connected); where $r(x)$ is the spectral
radius of $x$ in $A$.
 \item[(2)] $L(A)$ is locally arc connected.
\item[(3)] The singleton components of $L(A)$ are contained in the
centre of $A$ while the unboun\-d\-ed components are disjoint from
this centre.
 \item [(4)] For two distinct components $K_1$, $K_2$
of  $L(A)$; $\displaystyle d(K_1, K_2) \geq \rho(K_1, K_2) \geq
1$; where $\rho$ is the spectral distance. \item[(5)] The distance
of an unbounded component of $L(A)$ from the centre of $A$ is at
least $\frac{1}{2}$.
\end{enumerate}


\setcounter{paragraph}{0}
\section{Example}
 $A =
M_2(\raisebox{-0.5em}{$\stackrel{\bigtriangleup}{-}$})$, the
algebra of lower triangular $2\times 2$ matrices.
\[ L(A) = \left\{\left(\begin{array}{cc}
0 & 0 \\
0 & 0 \end{array} \right),  \left(\begin{array}{cc} 1 & 0 \\
0 & 1 \end{array} \right) \right\}\bigcup \left\{
\left(\begin{array}{cc} 1 & 0 \\
\alpha & 0  \end{array} \right) : \alpha \in K \right\}. \]
\[ R(A) = \left\{\left(\begin{array}{cc}
0 & 0 \\
0 & 0 \end{array} \right),  \left(\begin{array}{cc} 1 & 0 \\
0 & 1 \end{array} \right) \right\}\bigcup \left\{
\left(\begin{array}{cc} 0& 0 \\
\alpha & 1  \end{array} \right) : \alpha \in K \right\}. \]

\end{document}